\documentclass[11pt,reqno]{amsart}
\usepackage[a4paper, hmargin={2.7cm,2.7cm},vmargin={3.3cm,3.3cm}]{geometry}
\usepackage[english]{babel}
\usepackage{color}
\usepackage[noadjust]{cite}
\usepackage{amsmath}
\usepackage{amssymb}
\usepackage{amsfonts}
\usepackage{amsthm}
\usepackage{dsfont}
\usepackage{enumerate}
\usepackage{amsthm}
\usepackage{todonotes}
\usepackage{hyperref}
\usepackage{etoolbox, url}

\hypersetup{
	colorlinks   = true, 
	urlcolor     = blue, 
	linkcolor    = blue, 
	citecolor   = red 
}

\numberwithin{equation}{section}

\makeatletter
\@namedef{subjclassname@2020}{\textup{2020} Mathematics Subject Classification}
\makeatother

\makeatletter
\patchcmd{\ttlh@hang}{\parindent\z@}{\parindent\z@\leavevmode}{}{}
\patchcmd{\ttlh@hang}{\noindent}{}{}{}
\makeatother

\theoremstyle{plain}
\newtheorem{theorem}{Theorem}

\newtheorem*{conclusion}{Conclusion}

\theoremstyle{definition}

\theoremstyle{remark}

\title[Comment on `Comment on ....']{Comment on `Comment on `Affine density,
von Neumann dimension and a problem of Perelomov\textquotedblright }

\author{Lu\'{\i}s Daniel Abreu}

\address{Faculty of Mathematics, University of Vienna, Oskar-Morgenstern-Platz 1, A-1090 Vienna, Austria.}
\email{abreul22@univie.ac.at}

\begin{document}

\maketitle

\begin{abstract}
In the Letter `Comment on `Affine density, von Neumann dimension and a
problem of Perelomov' \url{https://doi.org/10.48550/arXiv.2211.04879}, by
Prof. J. L. Romero, it is claimed that the main theorem of Ref2 := [Adv.
Math. 407, Article ID 108564, 22 p. (2022)] is included in the prior
research survey Ref1 := [Expo. Math., 40(2), 265-301, 2022]. We provide some
background, and demonstrate with hard facts that \emph{this claim is false},
as well as other claims in the Letter which, in most cases, either result
from wrong mathematical reasoning or objectively differ from what is written
in the published material. The author might be interested in publicly
retracting and apologising for his comments.
\end{abstract}

\section{\textbf{Introduction}}

This Letter is the exercise of a right of answer to the comments in \cite%
{Letter}. In \cite{Letter}, it has been claimed that \emph{the main theorem
of Ref2 := [Adv. Math. 407, Article ID 108564, 22 p. (2022)] is included in
the prior research survey Ref1 := [Expo. Math., 40(2), 265-301, 2022]}. We
will demonstrate that \emph{this claim is false}, by providing evidence that
there is not enough information in \emph{[Expo. Math., 40(2), 265-301, 2022] 
}to conclude any of the original results in \emph{[Adv. Math. 407, Article
ID 108564, 22 p. (2022)]. }We will also comment on other misquote claims of 
\cite{Letter}, which in most cases either follow from errors in the logical
reasoning, or objectively \emph{do not match what is written in the
published material} and therefore seem to have been concluded from a
superficial reading of our paper. We have cited \emph{Ref1 := [Expo. Math.,
40(2), 265-301, 2022]} and properly aknowledged all information relevant to
our work.

The letter is organized as follows. In the next section we provide the
required background on the wavelet transform and define the spaces that are
fundamental for the discussion, as well as some literature review about
these spaces. Readers familiar with the topic can directly go to the third
section, where we specifically argue that the claim in the abstract of \cite%
{Letter} is false. We isolate this from the other comments, since it
represents a \emph{plagiarism claim} that will be \emph{categorically
refuted with hard facts}. In the fourth section, the remaining contents of
the Letter \cite{Letter} are discussed. Evidence is provided that the
author's further claims either follow from \emph{wrong mathematical reasoning%
}, or \emph{do not conform with what is written in the published material}.
The Letter is finished with further comments and with a suggestion for a
friendly solution of this bizarre situation.

\section{\textbf{Background on the wavelet transform}}

Consider the unitary representation of the affine group $G_{a}=\mathbb{R}%
\times \mathbb{R}^{+}\equiv \mathbb{C}^{+}$, acting on the Hardy space $%
H^{2}\left( \mathbb{C}^{+}\right) $

\begin{equation*}
\pi (z)\psi (t):=\frac{1}{\sqrt{y}}\psi \left( \frac{t-x}{y}\right) ,\quad
z=x+iy\in \mathbb{C}^{+}\text{.}
\end{equation*}%
\emph{The continuous wavelet transform} of a function $f$ with respect to a
wavelet\ $\psi $ is defined, for every $z=x+iy\in \mathbb{C}^{+}$, as 
\begin{equation}
W_{\psi }f(z)=\left\langle f,\pi ({z})\psi \right\rangle _{H^{2}\left( 
\mathbb{C}^{+}\right) }\text{.}  \label{wavelet}
\end{equation}%
Consider the Daubechies-Paul \cite{DauTS} mother wavelets chosen from the
family $\{\psi _{n}^{\alpha }\}_{n\in \mathbb{N}_{0}}$, $\alpha >0$, defined
in terms of the Laguerre polynomials as follows 
\begin{equation}
(\mathcal{F}\psi _{n}^{\alpha })(\xi ):=\xi ^{\frac{\alpha }{2}}e^{-\xi
}L_{n}^{\alpha }(2\xi ),\quad L_{n}^{\alpha }(t)=\frac{t^{-\alpha }e^{t}}{n!}%
\left( \frac{d}{dt}\right) ^{n}(e^{-t}t^{\alpha +n}),\text{ }\xi >0\text{.}
\label{mother}
\end{equation}%
We will write $d\mu ^{+}(z)=y^{-2}d\mu _{\mathbb{C}^{+}}(z)$, where $d\mu _{%
\mathbb{C}^{+}}(z)\ $is the Lebesgue measure in $\mathbb{C}^{+}$. The image
of $H^{2}\left( \mathbb{C}^{+}\right) $ under the map%
\begin{equation*}
W_{\psi _{n}^{\alpha }}:H^{2}\left( \mathbb{C}^{+}\right) \rightarrow L^{2}(%
\mathbb{C}^{+},d\mu ^{+})
\end{equation*}%
is a space with reproducing kernel $k_{\psi _{n}^{\alpha }}(z,w)=\frac{1}{%
C_{\psi _{n}^{\alpha }}}W_{\psi _{n}^{\alpha }}\psi _{n}^{\alpha }({w}%
^{-1}.z),$ (where $.$ stands for the product on $G_{a}$)$\quad $and$\quad
k_{\psi _{n}^{\alpha }}(z,z)=\frac{\Vert \psi _{n}^{\alpha }\Vert _{2}^{2}}{%
C_{\psi _{n}^{\alpha }}}$, where $C_{\psi _{n}^{\alpha }}=2\pi \left\Vert 
\mathcal{F}\psi _{n}^{\alpha }\right\Vert _{L^{2}(%
\mathbb{R}
^{+},t^{-1}dt)}^{2}$. \ We will denote this space by $\mathcal{W}_{\psi
_{n}^{\alpha }}(\mathbb{C}^{+})$.

*\textsl{For }$n=0$* the space $\mathcal{W}_{\psi _{0}^{\alpha }}(\mathbb{C}%
^{+})$ is related to \emph{Bergman spaces of analytic functions} \cite{HKZ}
(see, for instance, section 2.1 in \cite{Fractal}). As proven in \cite%
{AnalyticWavelet}, this is the only choice leading to spaces of analytic
functions. It can be shown that the only subspace of $L^{2}(\mathbb{C}%
^{+},d\mu ^{+})$ invariant under the classical projective representation $%
\tau ^{\alpha }$ of the group $\mathrm{PSL}(2,\mathbb{R})$ is the Bergman
space of analytic functions. The von Neumann dimension of the space $%
\mathcal{W}_{\psi _{n}^{\alpha }}(\mathbb{C}^{+})$ *\textsl{for }$n=0$* was
computed for the first time in the projective case by Radulescu \cite[%
Section 3]{Radulescu}, using an extension of the Schur-type orthogonality
formulas to irreducible projective representations with square integrable
coefficients. A recent approach without requiring Schur-type orthogonality
formulas (and thus neither irreducibility nor square integrability) was
presented by Jones \cite{Jones}.

*\textsl{For }$n>0$* the spaces $\mathcal{W}_{\psi _{n}^{\alpha }}(\mathbb{C}%
^{+})$ are related to \emph{spaces of non-analytic functions (real analytic
and polyanalytic \cite{AF} of order }$n$\emph{) }and for special choices of $%
\alpha $, to spaces of Maass forms \cite{Maass} (see \cite{Duke} for a
modern reference), or hyperbolic Landau Levels \cite{Comtet,Mouayn,AoP}.
There are very few papers dealing with the wavelet spaces $\mathcal{W}_{\psi
_{n}^{\alpha }}(\mathbb{C}^{+})$. Motivated by the connections with Maass
forms, physical models \cite{AoP}, DPP's \cite{Ensemble}, the large sieve 
\cite{DL} and polyanalytic function theory \cite{AF}, the author of this
Letter has been participating on a program to study several aspects of these
spaces *for $n>0$*, started in \cite{AoP} and recently developed in \cite%
{DL,Ensemble,Affine}. This is a new research direction in wavelet analysis,
and fundamental parts of it, like the computation of the matrix
representation coefficients, have only been completed very recently in \cite%
{DL}.

In our paper under discussion \cite{Affine}, \emph{a new projective
representation of }$PSL(2,\mathbb{R})$\emph{\ acting on }$\mathcal{W}_{\psi
_{n}^{\alpha }}(\mathbb{C}^{+})$\emph{\ (*for all non-negative integer }$n$*%
\emph{)} has been introduced. In Theorem 2.5 of \cite{Affine}, the wavelets
leading to the required $\mathrm{PSL}(2,\mathbb{R})$ invariance are shown,
under mild conditions, to be the functions (\ref{mother}). The invariance is
then used to compute the von Neumann dimension of the space *$\mathcal{W}%
_{\psi _{n}^{\alpha }}(\mathbb{C}^{+})$ for $n>0$*, using the method of
Jones \cite{Jones}. As an application, in \cite{Affine} we obtain two
original results about frames in the spaces $\mathcal{W}_{\psi _{n}^{\alpha
}}(\mathbb{C}^{+})$ for $n>0$ (Theorem 2.1 and Theorem 2.8) and four
Corollaries with variations and reflecting the connections between different
topics. The case $n=0$ of our Theorem 2.1 in \cite{Affine} has been
considered before in \cite{Jordi}. This is acknowledged in the comments
before Corollary 2.4 in \cite{Affine}.

\section{\textbf{The abstract claim}}

The following results are discussed in \cite{Letter}:

\textbf{A}-The statement below, concerning a locally compact, second
countable unimodular group $G$. A discrete series $\sigma $-representation
of $G$ is \emph{defined by the authors in \cite{Jordi} as an irreducible
projective representation with square integrable coefficients}, the same
setting of Definition 3.1 in \cite{Radulescu}.

\begin{theorem}[{\protect\cite[Theorem 7.4]{Jordi}}]
Let $\Gamma \subseteq G$ be a lattice and let $(\pi ,\mathcal{H}_{\pi })$ be
a discrete series $\sigma $-representation of $G$ of formal dimension $%
d_{\pi }>0$.

\begin{enumerate}
\item[(i)] If $\pi |_{\Gamma }$ admits a cyclic vector, then $vol(G/\Gamma
)d_{\pi }\leq 1$.

\noindent In particular, if $\pi |_{\Gamma }$ admits a frame vector, then $%
vol(G/\Gamma )d_{\pi }\leq 1$.

\item[(ii)] If $\pi |_{\Gamma }$ admits a Riesz vector, then $vol(G/\Gamma
)d_{\pi }\geq 1$.
\end{enumerate}
\end{theorem}

\bigskip

\textbf{B}-The following definition and result.

In \cite{Affine}, we define the following \emph{new} projective
representation $\tau _{n}^{\alpha }$ of the group $\mathrm{PSL}(2,\mathbb{R}%
) $ acting on the space $\mathcal{W}_{\psi _{n}^{\alpha }}=\mathcal{W}_{\psi
_{n}^{\alpha }}(\mathbb{C}^{+})$, which is the image of the Hardy space $%
H^{2}\left( \mathbb{C}^{+}\right) $ under the wavelet transform with the
special window function $\psi _{n}^{\alpha }\in H^{2}(\mathbb{R})$ defined
in (\ref{mother}). 
\begin{equation}
\tau _{n}^{\alpha }(m^{-1})F(z)=\bigg(\frac{|cz+d|}{cz+d}\bigg)^{2n+\alpha
+1}F(m\cdot z),\quad z\in \mathbb{C}^{+},\;m=%
\begin{pmatrix}
a & b \\ 
c & d%
\end{pmatrix}%
\in \mathrm{PSL}(2,\mathbb{R}).  \label{projective}
\end{equation}%
Below is our Theorem 2.1 from \cite{Affine}.

\begin{theorem}[{\protect\cite[Theorem 2.1]{Affine}}]
Let $\Gamma \subset \mathrm{PSL}(2,\mathbb{R})$ be a Fuchsian group with
fundamental domain $\Omega \subset \mathbb{C}^{+}$ (of finite volume) and $%
F\in \mathcal{W}_{\psi _{n}^{\alpha }}(\mathbb{C}^{+})$. If $\{\tau
_{n}^{\alpha }(\gamma )F\}_{\gamma \in \Gamma }$ is a frame for $\mathcal{W}%
_{\psi _{n}^{\alpha }}(\mathbb{C}^{+})$, then 
\begin{equation}
\left\vert \Omega \right\vert \leq \frac{C_{\psi _{n}^{\alpha }}}{\Vert \psi
_{n}^{\alpha }\Vert _{2}^{2}}=\frac{4\pi }{\alpha }\text{,}
\label{eq:thm-sam}
\end{equation}%
where $|\Omega |$ is calculated via the measure $\mu $.

If $\{\tau _{n}^{\alpha }(\gamma )F\}_{\gamma \in \Gamma }$ is a Riesz
sequence for $\mathcal{W}_{\psi _{n}^{\alpha }}(\mathbb{C}^{+})$, then 
\begin{equation}
\left\vert \Omega \right\vert \geq \frac{C_{\psi _{n}^{\alpha }}}{\Vert \psi
_{n}^{\alpha }\Vert _{2}^{2}}=\frac{4\pi }{\alpha }\text{.}
\label{eq:thm-int}
\end{equation}
\end{theorem}

\begin{itemize}
\item After stating these two results, the first comment of \cite{Letter} is:%
\begin{equation*}
\end{equation*}
\end{itemize}

`\emph{Let us start by mentioning that in \cite[Example 9.2]{Jordi} we
discuss at length how Theorem 1 applies to the so-called holomorphic
discrete series of }$PSL(2,\mathbb{R})$\emph{, which exhaust up to complex
conjugation and projective unitary equivalence all square integrable
irreducible representations of }$PSL(2,\mathbb{R})$\emph{\ (including (\ref%
{projective}))}.'%
\begin{equation*}
\end{equation*}

\begin{itemize}
\item This comment is \emph{mathematically wrong}. While the so-called
holomorphic discrete series of $\mathrm{PSL}(2,\mathbb{R})$ do exhaust up to
complex conjugation and projective unitary equivalence all \emph{square
integrable irreducible representations} of $\mathrm{PSL}(2,\mathbb{R})$,
this obviously does not include (\ref{projective}), which is only a\emph{\
projective representation }for \emph{non-integer} $\alpha $. See Bargmann 
\cite[5g, II]{Bar}.

\item For integer $\alpha $, the unitary equivalence between the projections
simplifies the application of Theorem 1, since the irreducibility and square
integrability (that we didn't need to use with our method-see Theorem A
below) of the representation follows automatically. But for general $\alpha $
this cannot be assured without a proof. See the very last point of this
Letter to see why this cannot be assured by general results for the
projective representation restricted to a Fuchsian group, not even by
lifting to the universal cover (the relevant cohomology obstruction only
vanishes for some Fuchsian groups). 
\begin{equation*}
\end{equation*}
\end{itemize}

Then, spread along the paper, there are some comments aiming at supporting
the following claim from the abstract:%
\begin{equation*}
\end{equation*}

`\emph{We point out that the main theorem of Ref2 := [Adv. Math. 407,
Article ID 108564, 22 p. (2022)] is included in the prior research survey
Ref1 := [Expo. Math., 40(2), 265-301, 2022]}.'%
\begin{equation*}
\end{equation*}

\begin{itemize}
\item By stating that the result is contained in the mentioned prior
research, the author is making a \emph{plagiarism claim}. We will
demonstrate that this claim is false, by showing that there is not enough
information in [Expo. Math., 40(2), 265-301, 2022] (nor, as far our
knowledge goes, in other prior research) to obtain Theorem 2. This requires
the following half-page argumentation.
\end{itemize}

To obtain Theorem 2 from Theorem 1 it is required to:

\begin{enumerate}
\item Find a \textsl{new} projective representation $\tau _{n}^{\alpha }$ of 
$\mathrm{PSL}(2,\mathbb{R})$ naturally acting on the spaces $\mathcal{W}%
_{\psi _{n}^{\alpha }}(\mathbb{C}^{+})$, because the classic projective
representation of $\mathrm{PSL}(2,\mathbb{R})$ on analytic Bergman spaces
does not leave $\mathcal{W}_{\psi _{n}^{\alpha }}(\mathbb{C}^{+})$ invariant
for $n>0$.

\item Show that the reproducing kernel space $\mathcal{W}_{\psi _{n}^{\alpha
}}(\mathbb{C}^{+})$ is invariant under $\tau _{n}^{\alpha }$.

\item Show that the representation $\tau _{n}^{\alpha }$ has square
integrable coefficients.

\item Show that the representation $\tau _{n}^{\alpha }$ is irreducible.

\item Compute the von Neumann dimension $\dim _{\pi |_{\Gamma }}\mathcal{W}%
_{\psi _{n}^{\alpha }}(\mathbb{C}^{+})\mathcal{=}vol(G/\Gamma )d_{\pi }$ of
the space *$\mathcal{W}_{\psi _{n}^{\alpha }}(\mathbb{C}^{+})$ for $n>0$*,
with the representation $\tau _{n}^{\alpha }$ (the $d_{\pi }=d_{\pi }(n)$
was first computed for $n>0$ by us in \cite{Affine}).
\end{enumerate}

The items (1)-(5) are \textsl{not contained in} [Expo. Math., 40(2),
265-301, 2022] (nor, as far our knowledge goes, in other published
research). It is actually enough to say that the projective representation $%
\tau _{n}^{\alpha }$ has been considered for the first time in our paper
[Adv. Math. 407, Article ID 108564, 22 p. (2022)]. Because without the
explicit form of the projective representation $\tau _{n}^{\alpha }$ it is
obviously \emph{not possible to state (1)-(5) (not even to define the
coefficients of the frame in Theorem 2)}, and therefore not possible to
prove these required properties. It should be noticed that the spaces $%
\mathcal{W}_{\psi _{n}^{\alpha }}(\mathbb{C}^{+})$ for $n>0$ are not even
mentioned in [Expo. Math., 40(2), 265-301, 2022].

Our approach in \cite{Affine} is simpler than the one using Theorem 1,
because it does not require to prove irreducibility and square
integrability. It uses two elementary results on von Neumann algebras,
which, according to Jones \cite[Theorem 3.5, (viii),(ix)]{Jones} are due to
Murray and von Neumann \cite{MurrayvonNeumann} (1936). These results are
stronger (and off course, older) than Theorem 1 i) and, combined with the
relation between Riesz and separating vectors (a novel result from \cite%
{Jordi} that we use and properly quote), imply Theorem 2 requiring only
conditions (1), (2) and (5). We invite the readers to our simple proof of
Theorem 2.1 in section 4.3 of \cite{Affine}. For a sample of the strength of
the results we used, observe that, for a general von Neumann algebra $M$
acting on a Hilbert space $\mathcal{H}$, the following holds:

\textbf{Theorem A [}Murray and von Neumann \cite{MurrayvonNeumann} (1936)]. $%
M$\emph{\ admits a cyclic vector if and only if} $\dim _{M}\mathcal{H}\leq 1$%
.

Take $M=\pi |_{\Gamma }$. Then $\dim _{M}\mathcal{H=}vol(G/\Gamma )d_{\pi }$
and that's it. We obtain, from a\emph{\ }86 years old result:

\textbf{Theorem A' [}Murray and von Neumann \cite{MurrayvonNeumann} (1936)]. 
$\pi |_{\Gamma }$\emph{\ admits a cyclic vector if and only if} $%
vol(G/\Gamma )d_{\pi }\leq 1$.

We observe that\emph{\ \textbf{Theorem A' }is stronger than Theorem 1 i): }%
it gives the same conclusions without the discrete series $\sigma $%
-representation  assumptions on the projective representation. \textsl{%
Therefore}, \textsl{in the final part of our argument, we use Theorem A
instead of Theorem 1 (this seems to be the main concern of the author of 
\cite{Letter}) because: }{\large it's simpler, stronger, and older.}

\begin{conclusion}
The plagiarism claim made in the abstract of \cite{Letter} is \textsl{false}.
\end{conclusion}

\section{\textbf{Comments on the Letter}}

\subsection{On section 1.1}

In this section the author presents a computation of $d_{\pi }$, alternative
to ours in \cite{Affine}. There are not many comments to be made, since this
was not included in [Expo. Math., 40(2), 265-301, 2022]. But perhaps one
should draw to the attention of the author that his proof is incomplete. The
author is confused about the definition of discrete series $\sigma $%
-representation from his own work with van Velthoven: after sorting the
Jargon and notation in the beginning of Section 1.1, this is defined simply
as a projective representation and the square integrability and
irreducibility, fundamental in the definition of discrete series $\sigma $%
-representation are forgotten. For this reason, the author didn't prove the
Schur-type orthogonality relations \cite[(1.5)]{Letter}, corresponding to
Definition 3.1 in \cite{Radulescu}, which only hold for\emph{\ irreducible }%
projective representations with \emph{square integrable coefficients }(in
Romero and van Velthoven's parlance, this is exactly a discrete series $%
\sigma $-representation). Neither the irreducibility nor the square
integrability of the representation $\tau _{n}^{\alpha }$ is shown in \cite%
{Letter}. This amounts to some extra work and is essential for the suggested
approach. We have actually done this work, but since our approach in \cite%
{Affine} does not require these two properties, the proof is omitted and
left as an exercise to the interested reader.

\subsection{On section 2}

In this section the author reproduces some of the proofs of \cite{Jordi},
perhaps to simplify the work of the reader. I have no comments to make about
this section.

\subsection{On section 3.2, `Our work'}

The section starts with a variation on a previous theme:

\begin{equation*}
\end{equation*}

`\emph{Firstly, our research survey \cite{Jordi} discusses at length the
application of Theorem 1 - and its counterpart concerning the existence of
coherent systems - to the holomorphic series of }$PSL(2,\mathbb{R})$\emph{,
which exhaust up to complex conjugation and projective unitary equivalence
all square integrable irreducible representations of }$PSL(2,\mathbb{R})$. 
\emph{The statements in Theorem 2 are of course invariant under unitary
equivalence.}'%
\begin{equation*}
\end{equation*}

\begin{itemize}
\item For those who missed the previous related comment in Section 3, here
it goes again. While the so-called holomorphic discrete series of $\mathrm{%
PSL}(2,\mathbb{R})$ do exhaust up to complex conjugation and projective
unitary equivalence all \emph{square integrable irreducible representations}
of $\mathrm{PSL}(2,\mathbb{R})$, this obviously does not include (\ref%
{projective}), which, for non-integer $\alpha $, is only a\emph{\ projective
representation. }Thus, the above comment has no relevance for the current
discussion. See also the last point of this Letter.%
\begin{equation*}
\end{equation*}
\end{itemize}

Then,%
\begin{equation*}
\end{equation*}

`\emph{the special case of Theorem 2 with }$F=W_{\psi _{n}^{\alpha }}\psi
_{n}^{\alpha }$\emph{\ is attributed to us in \cite{Affine}, whereas the
reader is led to understand that the generalization to arbitrary F is the
original contribution of [2]. While, indeed, in [\cite{Jordi}, Section
9.1.3, \textquotedblleft Perelomov's problem with respect to other special
vectors\textquotedblright ] we discuss the application of Theorem 1 to the
special vectors in question, the very same arguments apply to any other
generator}'%
\begin{equation*}
\end{equation*}

\begin{itemize}
\item This is a total mess. The original contribution of Theorem 2 is the
extension to general $n$, not to general $F$. What we attribute to \cite%
{Jordi} is the special case of Theorem 2 with $F=W_{\psi _{n}^{\alpha }}f$ \
(case $n=0$, general $F$! -this implies the result for $F=W_{\psi
_{n}^{\alpha }}\psi _{n}^{\alpha }$ - the general $n$ case $F=W_{\psi
_{n}^{\alpha }}f$ is not possible to obtain without the steps (1)-(5), in
particular needs the von Neumann dimension of $\mathcal{W}_{\psi
_{n}^{\alpha }}(\mathbb{C}^{+})$). Because what is proved in \cite{Jordi} is
the case $n=0$ of Theorem 2 (as recognized by us in the comments before
Corollary 2.4 in \cite{Affine}, which is the result we attributed to \cite%
{Jordi}). More precisely, it is shown in \cite[Section 9.1.2-9.1.3]{Jordi}
that if $\{\pi _{\alpha }\left( \gamma \right) g\}_{\gamma \in \Gamma }$ is
a frame (and therefore complete) in $A_{\alpha }^{2}\left( \mathbb{C}%
^{+}\right) $, where $A_{\alpha }^{2}\left( \mathbb{C}^{+}\right) $ stands
for the \emph{Bergman space of analytic functions }and $\pi _{\alpha }$ for
the classical projective representation of $\mathrm{PSL}(2,\mathbb{R})$ in $%
A_{\alpha }^{2}\left( \mathbb{C}^{+}\right) $, then the conclusion of
Theorem 2 holds. Since the Bergman space is a weighted version of $\mathcal{W%
}_{\psi _{0}^{\alpha }}(\mathbb{C}^{+})$ (see for instance section 2.1 in 
\cite{Fractal}), this is precisely the statement of Theorem 2 for $n=0$. The
author of \cite{Letter} should read more carefully the \textsl{two} papers
under discussion.%
\begin{equation*}
\end{equation*}
\end{itemize}

Then:

\bigskip

\emph{Thirdly, with a certain physical motivation, the authors of \cite{AoP}
consider the above-mentioned affine coherent states and conclude that the
density condition }%
\begin{equation*}
|\Omega |\leq \frac{4\pi (n+1)}{\alpha }
\end{equation*}%
\emph{is necessary for completeness.} \emph{In [\cite{Jordi}, Section 9.1.3]
we point out that coupling theory, as embodied in Theorem 1, offers the
sharper bound}%
\begin{equation}
|\Omega |\leq \frac{4\pi }{\alpha }  \label{bound}
\end{equation}

\begin{itemize}
\item This is \emph{true}! I have given before my congratulations to the
authors of \cite{Jordi} for improving our 2015 results in \cite{AoP},
despite the general case $\alpha $ not being completely executed (see the
last remarks). And we have quoted this result before our formulation of
Corolary 2.4 in \cite{Affine}: `\emph{As recently observed by Romero and
Velthoven in a different formulation...}'. We are actually quoting more than
what is proved in \cite{Jordi}, since we also attribute the Riesz sequence
condition. This was not proved in \cite{Jordi} but it follows from the
material in \cite{Jordi} (in contrast with our Theorems 2.1, 2.5 and 2.8 in 
\cite{Affine}, which, as already explained, do not; the estimate (\ref{bound}%
) is the case $n=0$ of our theorem 2.1; it was fully executed in \cite{Jordi}
for the case of $\alpha $ integer). In the paper with Balazs, de Gosson and
Mouayn \cite{AoP}, where the problem was considered for the first time, we
have used automorphic forms and this required introducing multiple zeros,
whence the factor $(n+1)$. The authors of \cite{Jordi} deserve all the
credit for this \emph{amazing} improvement of our early results (on our
behalf, \cite{AoP} was the first paper on the topic, and sometimes first
papers don't achieve best possible results).
\end{itemize}

\begin{equation*}
\end{equation*}%
Continuing,%
\begin{equation*}
\end{equation*}

\emph{The final remark in [\cite{Affine}, Section 2] on the
\textquotedblleft physical relevance of the results\textquotedblright\ may
give the impression that this observation is original to \cite{Affine}; we
invite the interested reader to look into [\cite{Jordi}, Section 9.1.3]}.'%
\begin{equation*}
\end{equation*}

\begin{itemize}
\item There is a\emph{\ scandalous mismatch} between this claim and what is
written in our paper. The author of \cite{Letter} is once again advised to
read more carefully the papers under discussion before posting public
comments of this nature. Our final remark on the `physical relevance of the
results' \cite[Page 9, line -6]{Affine} concludes as follows: `\emph{%
Corollary 2.4 provides estimates} for the number of particles distributed by
the model on such a region'. It turns out that \emph{Corollary 2.4 is
precisely the result we attributed to Romero and Velthoven }\cite{Jordi} and
this is the only result quoted in the whole paragraph! The reader may want
to know \emph{how} can this `\emph{give the impression that this observation
is original to }\cite{Affine}'.

\item We are actually using Corollary 2.7 from \cite{Affine}, which uses the
connection to the Maass forms operator, but since this amounts to a
specialization of the parameters in Romero and Velthoven's result, we were
extremely careful and, to respect their original contribution, we quoted
instead Corollary 2.4.
\end{itemize}

\section{\textbf{Further comments}}

\begin{itemize}
\item Our work was strongly influenced by Jones last paper \cite{Jones},
which offers new methods, simpler than those in the previous literature.
This was intended to be the first paper on what he envisioned as a new field
called \emph{Applied von Neumann algebra}. I leave the YouTube talk link for
those interested\emph{\ \url{https://www.youtube.com/watch?v=s7OvmKv4xfg}. }%
Before our paper, the relevant problems have only been considered in the
analytic Bergman setting, the natural environment of automorphic forms. We
offer in \cite{Affine} a contribution for this field, by extending the
Bergman space setting of \cite{Jones} to an infinite sequence of wavelet
spaces including the spaces of Maass forms and computing the von Neumann
dimension for all of them. These spaces are central in mathematics, since
the Maass Laplacian is the spinal cord of Selberg's trace formula (see the
clear presentation of Patterson \cite{Patterson}). Our projective
representation (\ref{projective}) didn't come from outer space, it was
designed to match the automorphic coefficient of Maass forms \cite%
{Duke,Maass}. This opens a series of questions, like constructing the
functions vanishing on the orbits of Fuchsian groups, whose existence is
assured by Theorem 2.1 (see Corollary 2.7 in \cite{Affine}). Such
construction is expected to involve Maass forms in a nontrivial way. Deep
problems arise if one attempts to extend other problems treated by Jones in 
\cite{Jones}, to our setting.

\item \emph{On section 3.1}. There are some comments about the status of
Perelomov problem. We have clarified this in our Remark 2.6 of \cite{Affine}%
. Abstract results can also be seen as a solution, this is really a personal
perspective, and the in \cite{Letter} quoted Bekka's claim from \cite{Bekka}
goes into this direction, but does not invalidate ours, which is concrete
and even classify the concrete wavelets for which the problem can be
considered (Theorem 2.5 in \cite{Affine}). As we point out in Remark 2.6 of 
\cite{Affine}, the really difficult part is still open. We invite the
interested reader to read the Remark, the fascinated reader to try to prove
what is still open, and the interested and fascinated reader to give a go at
Seip's conjecture mentioned in the introduction. I would love to see a proof
before leaving this existence.

\item \emph{On section 3.1}. \emph{Second paragraph. }I had to read our
whole paper to see where in `\emph{Sections 1.1 and 1.2 of \cite{Affine}
elaborate on the many technical challenges that the exponential growth of
balls in the hyperbolic metric'.} Besides the introductory sentence, and the
mention in Section 1.1 that Landau-type \cite{Landau,Mirko}, comparison \cite%
{HK} and analytic \cite{Seip} methods do not work here, nothing else is
mentioned. Obviously, the interest of these methods is to circunvent such
difficulties, which indeed do not show up in this approach. We have invested
efforts in giving a broad idea of the state of the art of the motivating
problems in the introduction.

\item The author of \cite{Letter} seems particularly fixed in our Theorem
2.1, suggesting it is not new (something already clarified in this Letter)
and calling it `main result'. We don't make such distinctions in our paper.
We offer three theorems (2.1, 2.5, 2.8) and six corollaries. My favorite is
Theorem 2.5, my coauthor prefers Theorem 2.8. All readers are invited to
read \cite{Affine} and make their own hit parade.

\item The constructive criticism of this Letter towards \cite{Letter} does
not include Romero and van Velthoven's work \cite{Jordi}. \cite{Jordi} is a
beautiful piece that goes beyond its survey intents, with original
contributions that were used and properly quoted in our work, and which
fully achieves the authors' stated goal in the introduction, since it really
`motivates the non-specialist to delve deeper into operator-algebraic
methods'. It also does a great job in providing a proof of Theorem 1 by
elementary methods inspired in `Janssen's classroom proof ' \cite{janssen}.
Therefore, we also invite all readers to read \cite{Jordi}.

\item Despite the praise in the previous paragraph, I need to point out a
small mistake in Example 9.3 of \cite{Jordi}, which seems to be one of the
many sources of\ the confusion that led the author of the Letter \cite%
{Letter} to write strange things, since the same mistake appears twice in 
\cite{Letter}, as already pointed out above. For non-integer $\alpha $, $\pi
_{\alpha }$ is just a projective representation, therefore not in the
holomorphic discrete series (as claimed in Example 9.3 of \cite{Jordi}),
since it has been shown by Bargmann \cite[5g, II]{Bar} that this implies $%
\alpha $ integer (see also \cite[3.3]{delaHarpe}). Thus, the orthogonality
relations cannot be obtained from the holomorphic discrete series. Thus, the
argument of the authors' only assures the bound (\ref{bound}) for $\alpha =n$
(thus, to be fully rigorous, only improve our results in \cite{AoP} in these
special cases). The general case of non-integer $\alpha $ is not completely
executed. We have ignored this, and attributed the result with general $%
\alpha $ to \cite{Jordi} because, fortunately, this is not a serious
mistake: it can be corrected by replacing the concept of holomorphic
discrete series by irreducible square integrable projective representations
with square integrable representations (Definition 3.1 in \cite{Radulescu},
this is what the authors call discrete series $\sigma $-representation), and
then using the results of Section 3 of \cite{Radulescu}, the claimed
Schur-type orthogonality holds, and the conclusion of the authors' in
Example 9.3 of \cite{Jordi} follows.

\item For a crash overview of the holomorphic discrete series in a signal
analysis context, I would recommend Feichtinger-Gr\"{o}chenig \emph{coorbit
zero} \cite[Example 7.3]{FG}, where the case of discrete series ($\alpha =n$%
) is worked out and, after this, the possibility of extension to projective
representations is mentioned.

\item I invite the author of the Letter \cite{Letter}, after reading more
carefully the material involved, for a\emph{\ }\textsl{blackboard
mathematical discussion} about the topics involved in these papers and
Letters, at NuHAG's discussion room at the University of Vienna, together
with two Senior scientists from our group with the capacity of following the
contents and arbitrate between facts and viewpoints in case this is needed.
I believe Prof. Feichtinger and Prof. Gr\"{o}chenig are more than qualified
to assist us in this task. We can later invite Michael Speckbacher and Jordi
van Velthoven for the discussion (and also present the material in NuHAG's
weekly seminar) but, at this stage, it would be ideal to first clarify
things between the two of us.

\item The author might be interested in retracting the comments of \cite%
{Letter} and to voluntarily making a public apology to the authors of \cite%
{Affine} in arXiv, under the title \emph{Retraction and apology for }`\emph{%
Comment on `Affine density, von Neumann dimension and a problem of Perelomov}%
\textquotedblright .

\item A further remark about the projective representations. For groups with
trivial Lie algebra cohomology as the case of the whole $SL(2,\mathbb{R})$,
one could think about using Bargmann's theorem \cite{Bar1} to lift the
projective representations to representations in the universal cover, where
they were classified by Pukanszky \cite{Puk} and Sally \cite{Sally}. But
this is \emph{not possible} to do for the \emph{restriction of the
projective representation to an arbitrary Fuchsian group}. Indeed, in
Theorem A1 in \cite[Appendix A]{Jones}, Jones proved that, for a Fuchsian
group associated with a Riemann surface of genus $g$, the cohomology
obstruction only vanishes (a requirement to lift the representation) if $%
\alpha $ is an integer multiple of $\frac{1}{g-1}$. Thus, for a Fuchsian
group of genus $g=2$, this requires $\alpha $ to be an integer and we are
back in the holomorphic discrete series. Thus, there is no possibility of
assuring unitary equivalences between the results of Theorem 2.1 by general
results. The action of the projective representation (\ref{projective}) must
be considered in its explicit form, acting on each Fuchsian group, in its
associated space $\mathcal{W}_{\psi _{n}^{\alpha }}(\mathbb{C}^{+})$. Thus,
from any mathematical perspective, (\ref{projective}) is a \emph{new
projective representation of }$PSL(2,R)$.
\end{itemize}

Finally, it is with a profound sadness that I am forced to use a scientific
repository to write a Letter like this. But the unannounced posting of \cite%
{Letter} left no other choice. A plagiarism allegation from someone with the
credibility of Prof. Romero, could easily raise doubts among my colleagues
and end up with my career and reputation, as well of my coauthor's, if left
unanswered and unclarified. The situation is so unusual in mathematics that
I doubt of the existence of a system prepared to solve it. The only possible
definitive solution is to appeal to a consensus building among experts in
mathematics.%
\begin{equation*}
\end{equation*}

\textbf{Aknowledgements. }\emph{The Letter has been written at the author's
sole responsibility. We thank the author of \cite{Letter} for the correction
of a typo in the definition of the representation (\ref{projective}) and a
small mistake in the statement of Theorem 2.1 (as the author correctly
pointed out, the compactness condition is not necessary). We will include
these corrections, and a few more in a corrigendum to be submitted (the
statements in these Letter already include them, namely a correction of the
Nyquist rate, }$\frac{2}{\alpha }$\emph{\ should be replaced by }$\frac{4\pi 
}{\alpha }$\emph{\ an unfortunate mistake inherited from not taking into
account the normalization of the Fourier transform used in \cite{Affine}).
The light notes of irony in one or two points of this letter are intended to
soften down the tone of \cite{Letter}, which conveys a certain loftiness.
Furthermore, the subject matter at stake is too serious to be taken
seriously.}

\end{document}